\input epsf
\documentstyle{amsppt}
\pagewidth{6.48truein}
\pageheight{9.6truein}\vcorrection{-0.6in}
\TagsOnRight
\NoRunningHeads
\catcode`\@=11
\def\logo@{}
\footline={\ifnum\pageno>1 \hfil\folio\hfil\else\hfil\fi}
\topmatter
\title A short conceptual proof of Narayana's path-counting formula
\endtitle
\vskip-0.2in
\author Mihai Ciucu\endauthor
\thanks Research supported in part by NSF grant DMS-1101670 and DMS-1501052.
\endthanks
\affil
  Department of Mathematics, Indiana University,
  Bloomington, Indiana 47405
\endaffil
\abstract
We deduce Narayana's formula for the number of lattice paths that fit in a Young diagram as a direct consequence of the Gessel-Viennot theorem on non-intersecting lattice paths.
\endabstract
\endtopmatter

\document

\def\mysec#1{\bigskip\centerline{\bf #1}\message{ * }\nopagebreak\bigskip\par}

\def\myref#1{\item"{[{\bf #1}]}"}

\def\pf{{\it Proof.\ }}

\def\epf{\hfill{$\square$}\smallpagebreak}

\def\cite#1{\relaxnext@
  \def\nextiii@##1,##2\end@{[{\bf##1},\,##2]}%
  \in@,{#1}\ifin@\def\next{\nextiii@#1\end@}\else
  \def\next{[{\bf#1}]}\fi\next}
\def\proclaimheadfont@{\smc}

\def\pf{{\it Proof.\ }}

\define\Z{{\Bbb Z}}


\define\twoline#1#2{\line{\hfill{\smc #1}\hfill{\smc #2}\hfill}}
\define\twolinetwo#1#2{\line{{\smc #1}\hfill{\smc #2}}}
\define\twolinethree#1#2{\line{\phantom{poco}{\smc #1}\hfill{\smc #2}\phantom{poco}}}
\define\threeline#1#2#3{\line{\hfill{\smc #1}\hfill{\smc #2}\hfill{\smc #3}\hfill}}

\def\mypic#1{\epsffile{#1}}



\define\DT{1}
\define\GV{2}
\define\Kr{3}
\define\Na{4}
\define\Natwo{5}




\vskip-0.1in

Let $\lambda$ and $\mu$ be two partitions so that $\mu\subset\lambda$, and consider the skew Young diagram $\lambda/\mu$ (see Figure 1 for an example). Let $N(\lambda/\mu)$ be the number of minimal lattice paths on $\Z^2$ contained in this skew Young diagram from its southwestern to its northeastern corner. We give a short conceptual proof for the following extension of Narayana's path-counting formula \cite{\Na} due to  Kreweras \cite{\Kr} (see \ also \cite{\Natwo, Ch.\,II}; the special case $\mu=\emptyset$ yields Narayana's formula). 

\proclaim{Theorem 1 (Kreweras \cite{\Kr})} Let $n$ be the number of parts of $\lambda/\mu$. Then
$$
N(\lambda/\mu)
=
\det\left({\lambda_j-\mu_i+1\choose j-i+1}_{1\leq i,j\leq n}\right).\tag1
$$

\endproclaim

\pf The key to our proof is to consider the region $R$ on the triangular lattice corresponding to $\lambda/\mu$ indicated by the the outside contour in Figure 2 --- it is obtained from the Young diagram of $\lambda/\mu$ by affinely deforming it so that its unit squares become unit rhombi on the triangular lattice, and then translating the southeastern boundary one unit in the $-\pi/3$ polar direction. We argue that the number of tilings of $R$ by unit rhombi (a.k.a. lozenge tilings) is equal to both $N(\lambda/\mu)$ and the determinant in (1).

\midinsert
\twoline{\mypic{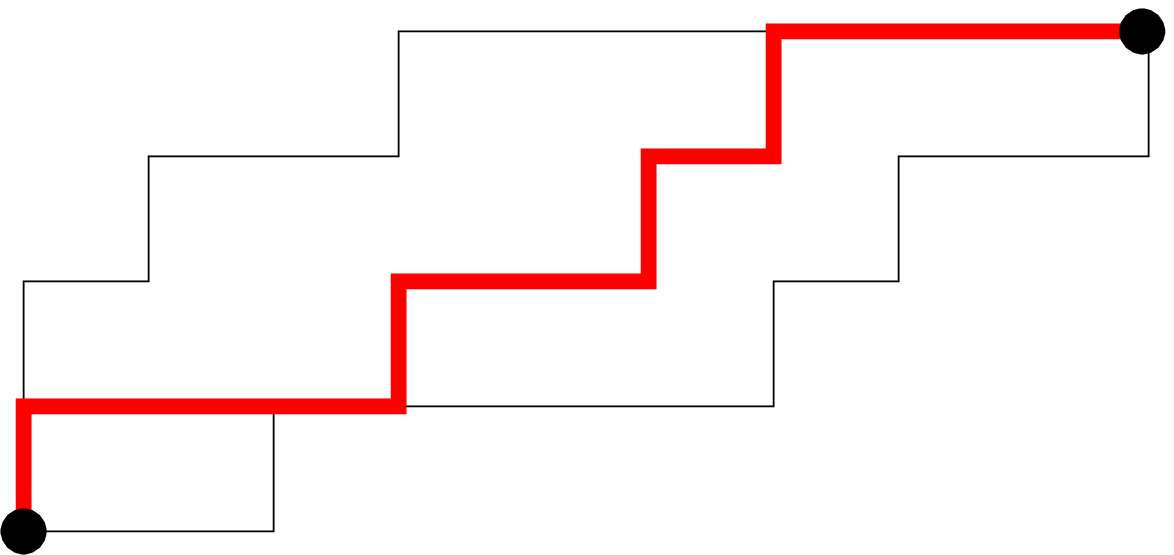}}{\mypic{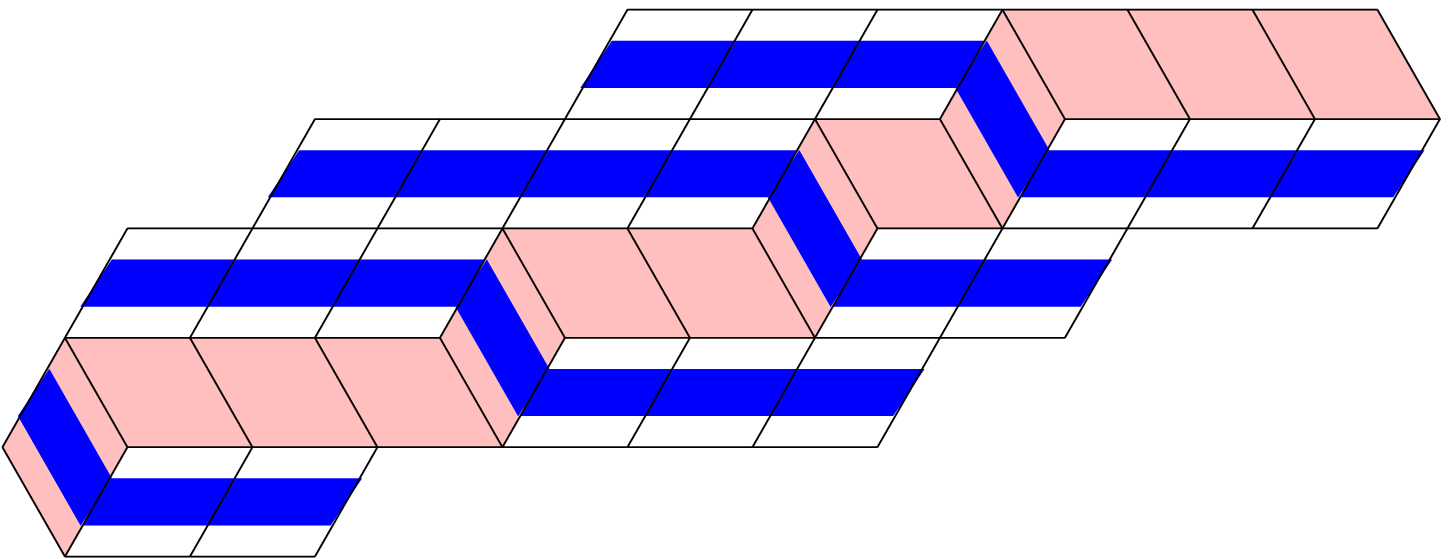}}
\twoline{Figure 1. {\rm The skew shape $(9,7,6,2)/(3,1)$.}}{Figure 2. {\rm The corresponding region $R$.}}
\endinsert

Indeed, recall that lozenge tilings of regions on the triangular lattice are in one-to-one correspondence with families on non-intersecting paths of rhombi (see \cite{\DT}). The latter can be chosen in three different ways, depending on whether the segments where the paths of lozenges start and end point in the $-\pi/3$, $\pi/3$ or $-\pi$ polar directions. For the region $R$, the first of these three ways yields a single path of rhombi (lightly shaded in Figure 2), which can be regarded as a lattice path in $\lambda/\mu$ connecting the southwestern and northeastern corners. On the other hand, the second way yields a family of $n$ non-intersecting paths of rhombi (shaded dark in Figure 2), which can be regarded as non-intersecting lattice paths on $\Z^2$. By the Gessel-Viennot theorem \cite{\GV}, their number is the determinant of the $n\times n$ matrix whose $(i,j)$-entry is the number of minimal lattice paths on $\Z^2$ from the $i$-th starting point to the $j$-th ending point (both counted from top to bottom). One readily checks that this is precisely the $(i,j)$-entry of the  matrix in (1). 
\epf

\mysec{References}
{\openup 1\jot \frenchspacing\raggedbottom
\roster

\myref{\DT}
  G. David and C. Tomei, The problem of the calissons,
{\it Amer\. Math\. Monthly} {\bf 96} (1989), 429--431.

\myref{\GV}
    I. M. Gessel and X. Viennot, Binomial determinants, paths, and hook length formulae,
{\it Adv. in Math.} {\bf 58} (1985), 300--321.

\myref{\Kr}
  G. Kreweras, Sur une classe de probl\`emes d\'enombrement li\'es au trellis des partitions de entiers, {\it Cahiers du Bur. Univ. de Rech. Op\'er.} {\bf 6} (1965), 5--105. 

\myref{\Na}
  T. V. Narayana, A combinatorial problem and its application to probability theory, {\it J. Indian Soc. Agr. Stat.} {\bf 7} (1955), 169--178.

\myref{\Natwo}
  T. V. Narayana, ``Lattice path combinatorics with statistical applications,'' {\it Mathematical Expositions No. 23}, University of Toronto Press, Toronto, 1979.

\endroster\par}

\enddocument